\documentstyle[12pt]{article}
\newlength{\abstractwidth}
\renewcommand{\thefootnote}{\fnsymbol{footnote}}
\renewcommand{\thanks}[1]{\footnote{#1}} 
\newcommand{\starttext}{ \setcounter{footnote}{0}
\renewcommand{\thefootnote}{\arabic{footnote}}}

\newcommand{\be}{\begin{equation}}
\newcommand{\bea}{\begin{eqnarray}}
\newcommand{\eea}{\end{eqnarray}} \newcommand{\ee}{\end{equation}}
 
 \def\ba{\begin{eqnarray}}
\def\ea{\end{eqnarray}}

\def\o{\omega}

\def\det{{\rm det}}

\def\log{\,{\rm log}\,}
\def\exp{\,{\rm exp}\,}

\def\o{\omega}

\def\o{\omega}

\def\ge{\geq}
\def\le{\leq}

\def\p{\partial}

\def\[{{\bf [}}
\def\]{{\bf ]}}

\def\ddbar{i\p\bar\p}

\def\mathbb{\bf}

\def\eqref{\ref}

\begin{document}
\starttext \baselineskip=18pt \setcounter{footnote}{0}
\newtheorem{theorem}{Theorem}
\newtheorem{lemma}{Lemma}
\newtheorem{corollary}{Corollary}
\newtheorem{definition}{Definition}
\newtheorem{conjecture}{Conjecture}
\newtheorem{proposition}{Proposition}


\begin{center}
{\Large \bf STABILITY ESTIMATES FOR THE COMPLEX MONGE-AMP\`ERE AND HESSIAN EQUATIONS
\footnote{Work supported in part by the National Science Foundation under grant DMS-1855947.}}

\medskip
\centerline{Bin Guo, Duong H. Phong, and Freid Tong}

\medskip

\begin{abstract}

{\footnotesize A new proof for stability estimates for the complex Monge-Amp\`ere and Hessian equations is given, which does not require pluripotential theory. A major advantage is that the resulting stability estimates are then uniform under general degenerations of the background metric in the case of the Monge-Amp\`ere equation, and under degenerations to a big class in the case of Hessian equations.
}

\end{abstract}

\end{center}

\baselineskip=15pt
\setcounter{equation}{0}
\setcounter{footnote}{0}

\section{Introduction}
\setcounter{equation}{0}

Stability estimates for a non-linear partial differential equation are estimates for how much the solution can vary, given the size of the variation of the right hand side.
Clearly, they are of great theoretical as well as practical importance. Such estimates had been obtained by Kolodziej \cite{K} for the complex Monge-Amp\`ere equation, by Dinew and Kolodziej \cite{DK} for complex Hessian equations, and by Dinew and Zhang \cite{DZ} for Monge-Amp\`ere equations when the background metric is not necessarily K\"ahler, but just big. In all cases, the proofs made extensive use of pluripotential theory and the background metric was fixed. It remained an open question whether these estimates can be established without pluripotential theory, and whether they can be made uniform under degenerations of the background metric, a situation which arises frequently in geometric applications.

\smallskip
In \cite{GPT}, the authors developed a method for obtaining sharp $L^\infty$ estimates for the complex Monge-Amp\`ere equation without pluripotential theory. As explained in greater detail there, the method of \cite{GPT} built on works of Wang, Wang, Zhou \cite{WWZ} and of Chen and Cheng \cite{CC}, particularly on the last two authors' idea of considering an associated complex Monge-Amp\`ere equation. It achieved the stated goal of giving an alternate PDE proof of $L^\infty$ estimates for the complex Monge-Amp\`ere equation, but it also went considerably beyond in, on one hand, applying to more general fully non-linear equations, and on the other hand, allowing the background metrics to degenerate. It can also give sharp gradient estimates \cite{GPT1}, improving on the estimates in e.g. \cite{PS, B, GL, Gp, H}.

\smallskip
The main goal of the present paper is to obtain stability estimates for the complex Monge-Amp\`ere and Hessian equations which are uniform under degenerations. We use the method of \cite{GPT}.
We recover in the process the stability estimates of \cite{K,DK, DZ}, this time without pluripotential theory. Our estimates are also uniform under general degenerations of the background metric in the case of the Monge-Amp\`ere equation, and under degenerations to a big class in the case of Hessian equations. Thus we answer in the positive both questions asked above.

\section{Statement of the main results}
\setcounter{equation}{0}

Let $(X,\omega)$ be a compact K\"ahler manifold, $\chi$ a closed and non-negative $(1,1)$-form, and set
\bea
\o_t=\chi+t\o,
\qquad t\in (0,1].
\eea
Let $f,h \in C^\infty$ are smooth  functions normalized by
$$\int_X e^f \omega^n = \int_X e^h \omega^n = \int_X \omega^n = 1,$$
and consider the following complex Hessian equations
\begin{equation}
\label{eqn:main}(\omega_t+ \ddbar u_t)^k\wedge \omega^{n-k} = c_t e^f \omega^n,\quad (\o_t+\ddbar v_t)^k \wedge \omega^{n-k} = c_t e^h \omega^n
\end{equation}
with the constants $c_t$ given by $c_t = \int_X \omega_t^k\wedge \omega^{n-k}$. We normalize $u_t$ and $v_t$ so that
$$\max_X (u_t - v_t) = \max_X(v_t - u_t).$$

\medskip
We will consider three cases: 
\bea
&{\mathrm  {\mathbf I}}:  &\nonumber k = n, \mbox{ and }t\in (0,1]\\
& {\mathrm {\mathbf {II}}}: &\nonumber  1\le k<n, \mbox{ and }\chi = 0, \, t = 1.\\
&{\mathrm {\mathbf {III}}}: &\nonumber 1\le k<n, \mbox{ and  } t\in (0,1],\, \chi \mbox{ is big i.e. $\int_X \chi^n >0$}
\eea
which correspond respectively to the Monge-Amp\`ere equations with degenerations, the $k$-th Hessian equation with a fixed background metric, and the $k$-th Hessian equation with degenerations. Different cases correspond to different choices of test functions, and constants. So we will treat the cases separately, when necessary.

\medskip
For each case, we will make the following assumptions and choice of constants, 
     \bea& {\mathrm  {\mathbf I}}: &\nonumber  \| e^{h}\|_{ L^1(\log L)^{p_1} (\omega^n)}, \| e^{f}\|_{ L^1(\log L)^{p_1} (\omega^n)} \le K,\, p_1> n\\
& {\mathrm {\mathbf {II}}}: &\nonumber  \| e^h\|_{L^{p_2}(\omega^n)} , \| e^f\|_{L^{p_2}(\omega^n)} \le K,\, p_2> \frac n k, \mbox{ and }q_2 = \frac{p_2 - 1}{p_2 -  n/ k} \\
& {\mathrm {\mathbf {III}}}:  &\nonumber  \| e^h\|_{L^{p_3}(\omega^n)} , \| e^f\|_{L^{p_3}(\omega^n)} \le K,\, p_3> \frac {n^2} k, \mbox{ and }q_3 = \frac{p_3 - 1}{p_3 -  n /k} 
      \eea
      for a fixed constant $K>0$.

The equations (\eqref{eqn:main}) admit smooth solutions by \cite{Y, DK1}. By \cite{GPT} (also \cite{K, DK}),  under the above assumptions, the oscillations ${\mathrm{osc}}u_t$ and ${\mathrm {osc}}v_t$ of the solutions are uniformly bounded independently of $t$ in all three cases ${\mathrm  {\mathbf I}}, {\mathrm  {\mathbf II}}$, and ${\mathrm  {\mathbf III}}$. 
Let  $\beta_0>1$ denote such an upper bound depending only on $n, k, \omega,\chi, p_a$ and $10K$. 

\medskip

To start with, we will define a positive function $\gamma_a(r)$ with $\gamma_a(r)\to 0$ as $r\to 0$. Each case has different choice of such a function.   
We define $\gamma_a$ case by case:
     \bea
     \label{cases}
     & {\mathrm  {\mathbf I}}: &\nonumber  \gamma_1(r) = r^{ \frac 1{\delta_1} - n + 1  },\mbox{ where } \delta_1 = \frac{p_1-n}{p_1n} < \frac 1  n\\
& {\mathrm {\mathbf {II}}}: &\nonumber  \gamma_2(r) = r^{(n+1)q_2 - n}\\ 
& {\mathrm {\mathbf {III}}}:  &\nonumber  \gamma_3(r) = r^{(n+1)q_3 - n}
      \eea

\begin{theorem}
\label{thm:main}
Let the assumptions and notations be as above. Then we have in all three cases listed in (\ref{cases})
\begin{equation}
\sup_{X} | u_t - v_t  |\le C \| e^f - e^h\|_{L^1}^{1/(n+3 + \sigma_a)},
\end{equation}
where in each case, $\sigma_a>0$ is the power of $r$ in $\gamma_a(r)$, i.e. $\gamma_a(r) = r^{\sigma_a}$, and $C$ is a constant depending only on  $n,k, \omega, \chi $, $K>0$ and $p_a$.

\end{theorem}

We observe that this theorem improves on all results known so far. More specifically in case {\bf I}, we get uniform stability for a degenerating family, and it does not even matter whether $\chi$ is big or not. Kolodziej \cite{K} proved this case for a fixed K\"ahler metric, and Dinew and Zhang \cite{DZ} proved it for a fixed big class. In case {\bf II},  we slightly sharpen the known stability result in \cite{DK}, where the RHS in the inequality is $\| e^f  - e^h\|_{L^{q'}}$ for some $q'>1$, while we are able to prove the inequality for $q' = 1$.  
In case {\bf III}, we obtain a uniform stability theorem for Hessian equations when the class remains big. This is completely new, and relies in particular on the uniform $L^\infty$ estimate in \cite{GPT} for solutions with  degenerating big classes. 

\smallskip

We note that the exponent $1/(n+3+\sigma_a)$ is not sharp in general. The sharp exponent can be obtained by replacing Lemma \ref{lemma Kol} below by a result from \cite{DZ}. We leave the details to the interested readers.

\section{Proof of Theorem \ref{thm:main}}
\setcounter{equation}{0}

For the proof, it is convenient to restate the theorem as follows. Under the above assumptions, if in each case we have
\begin{equation}\label{eqn:assumption basic}\| e^f - e^h\|_{L^1(\omega^n)}\le \gamma_a(r) r^{n+3},
\end{equation}
 then there is a small $r_0>0$ such that for all $0<r\le r_0$
$$\sup_X |u_t - v_t|\le C r,$$ for all $t$ listed in each case in (\ref{cases}) and $C$ depends only on $n,k, \omega, \chi $, $K>0$ and the corresponding $p_a$ in each case.
This is the version which we shall prove.

\bigskip

We begin with a lemma due to Kolodziej \cite{K}. The proof is almost identical to that in \cite{K}, but since we would like to avoid the use of pluripotential theory, some additional smoothing is needed in the proof, and we provide a full proof of this lemma.

\smallskip

Choose a small $\bar r_0\in (0, \frac 1 {10})$ such that $\gamma_a(\bar r_0) \bar r_0^{n} \le \frac 1 5$. We then fix an $0<r<\bar r_0$. We remark that all relevant constants are independent of $r$. Later on we will choose an even smaller $r_0>0$.

By switching the roles of $u_t$ and $v_t$ if necessary,  we may assume 
$$\int_{\{v_t\le u_t\}} (e^f+ e^h) \omega^n \le 1.$$

Denote $E_j := \{v_t \le u_t - j\beta_0 r\}$. The next lemma states that over the set $E_2$, the integral of $e^h$ is small.

\begin{lemma}\label{lemma Kol}
In each case $a = $ {\bf I, II, III}, we have 
$$\int_{E_2} e^h \omega^n\le C_0 \gamma_a(r) r^{n},$$ for some constant $C_0 = 1+ \frac{2}{ (\frac{3}{2})^{1/k} - 1}$.
 \end{lemma}
 
\noindent{\em Proof.}
We calculate
\begin{equation}\label{eqn:Kol 1}\int_{E_0 } e^h \omega^n = \frac 1 2 \int_{E_0} (e^f+e^h) + (e^h - e^f) \omega^n \le \frac 1 2 (1+ \frac 15) = \frac 3 5.
\end{equation}
Take a sequence of positive smooth functions  $\tau_j$ that converge uniformly to $\chi_{E_0}$ such that $\tau_j \equiv 1$ on $E_0$. Consider a sequence of smooth positive functions
$$e^{h_j} = \frac 3 2\tau_j e^h + c_j (1-\tau_j) e^h $$
where $c_j>0$ are chosen so that $\int_X e^{h_j}\omega^n = 1$. It is not hard to see from (\ref{eqn:Kol 1}) that for $j>>1$, $\frac 1{20} \le c_j \le 3.$ Hence  when $j>>1$

\medskip

$\bullet$ in case {\bf I}, we have  $\| e^{h_j}\|_{L^{1}(\log L)^{p_1}}\le 5K$,

\smallskip

$\bullet$ in case {\bf II}, we have  $\| e^{h_j}\|_{L^{p_2}}\le 5K$,

\smallskip

$\bullet$ in case {\bf III}, we have $\| e^{h_j}\|_{L^{p_3}} \le 5K$. 

\medskip

\newcommand{\bk}[1]{\Big(#1 \Big)}
\newcommand{\xk}[1]{\big(#1 \big)}

We solve the following Hessian equations which admit known to admit unique smooth solutions \cite{Y, DK1},
$$(\omega_t + \ddbar \rho_j)^k\wedge \omega^{n-k} = c_t e^{h_j} \omega^n,\quad \sup_X \rho_j = 0 \mbox{ and } \omega_t+\ddbar \rho_j\in \Gamma_k,$$ where $\Gamma_k$ is the usual open  convex cone in $k$-th Hessian equations.
By the choice of $\beta_0$, we have $-\beta_0 \le \rho_j \le 0$ (see \cite{GPT}). The following Newton inequality holds pointwise for any $1\le l \le k$
$$\omega_{t, u_t}^l \wedge \omega_{t, \rho_j}^{k-l}\wedge \omega^{n-k}\ge \bk{\frac{\omega_{t, u_t}^k\wedge \omega^{n-k}}{\omega^n}}^{l/k} \bk{\frac{\omega^k_{t, \rho_j} \wedge \omega^{n-k}}{\omega^n}}^{(k-l)/k} \omega^n.$$
Then on the set $E_0\backslash G $ where $G = \{e^f \le (1-r^2) e^h\}$, we have
$$\omega_{t, u_t}^l \wedge \omega_{t, \rho_j}^{k-l} \wedge \omega^{n-k}\ge c_t (1-r^2)^{l/k} (\frac 3 2)^{(k-l)/k} e^h \omega^n.  $$
It follows that on $E_0\backslash G$
     \bea\label{eqn:new one}
    \omega_{t,r\rho_j + (1-r) u_t}^k\wedge \omega^{n-k} & \nonumber =  &\nonumber \sum_{l=0}^k \frac{k!}{l!(k-l)!} r^{k-l} (1-r)^l \omega_{t, u_t}^l \wedge \omega_{t,\rho_j}^{k-l}\wedge \omega^{n-k} \\
   &\ge& \nonumber c_t \sum_{l=0}^k \frac{k!}{l!(k-l)!} r^{k - l} (1-r)^l (1-r^2)^{l/k} (\frac 3 2)^{(k-l)/k} e^h \omega^n\\
   & =  &c_t\xk{ (1-r)(1-r^2)^{1/k} + r (\frac  3 2)^{1/k}  }^k e^h \omega^n \ge  c_ t(1+b_0 r ) e^h \omega^n 
    \eea
where $b_0 = \frac 1 2 ((\frac{3}{2})^{1/k} - 1)>0$, since $r$ is chosen to be small. 

\medskip

Note that by (\ref{eqn:assumption basic})

$$r^2 \int_G e^h\omega^n \le \int_G (e^h - e^f)\omega^n \le \gamma_a(r) r^{n+3},$$ which implies \begin{equation}\label{eqn:G estimate}\int_G e^h\omega^n \le \gamma_a(r) r^{n+1}.\end{equation}

Adding the same constant to $u_t$ and $v_t$, we may assume without loss of generality $-\beta_0 \le u_t\le 0$. The following inclusion relation holds from the definition 
$$E_2 \subset E:=\{ v_t \le r \rho_j + (1-r) u_t - \beta_0 r  \}\subset  E_0.$$
All functions involved are smooth so by the comparison principle and (\ref{eqn:new one}),
    \bea
    c_t (1+b_0 r) \int_{E \backslash G} e^h \omega^n & \le &\nonumber  \int_E \omega_{t, r \rho_j + (1-r ) u_t}^k\wedge \omega^{n-k} \\
& \le &\nonumber  \int_E \omega_{v_t}^k \wedge \omega^{n-k}  = c_t \int_{E \backslash G} e^h \omega^n + c_t \int_G e^h \omega^n
     \eea
Combined with (\ref{eqn:G estimate}) this implies
$$\int_{E\backslash G} e^h \omega^n \le \frac{1}{b_0} \gamma_a(r) r^{n}$$
It follows that
$$\int_{E_2} e^h \omega^n \le \int_{E\backslash G} e^h \omega^n  + \int_G e^h \omega^n \le (1+ \frac 1 {b_0}) \gamma_a(r) r^n.$$
The Lemma is proved.

\bigskip
We now come to the proof of Theorem 1 proper. 
We  normalize $u_t$ as in the proof of Lemma \ref{lemma Kol}. For $s\ge 0$, we set $\Omega_s = \{v_t \le (1-r) u_t - 3 \beta_0 r - s\}$. Note that $\Omega_s\subset E_2$ for any $s\ge 0$.

\medskip

\medskip
We follow the same strategy as in \cite{GPT}.
We choose a sequence of smooth positive functions $\eta_j: {{\mathbb R}}\to{{{\mathbb R}}}_+$ such that 
\begin{equation}\nonumber
\eta_j(x) = x + \frac 1 j,\quad \mbox{ when }x\ge 0,
\end{equation}
 and $$ \eta_j(x) = \frac 1 {2j},\quad  \mbox { when } x\le -\frac {1}{j},$$ and $\eta_j(x)$ lies between $1/2j$ and $1/j$ for $x\in [-1/j, 0]$.  Clearly $\eta_j\to \eta_\infty(x) = x\cdot \chi_{{\mathbb R}_+}(x)$ pointwise as $j\to \infty$.

 \medskip
 
 We solve the complex Monge-Amp\`ere equations
             $$(\omega_t+ \ddbar \psi_j)^n = c_t^{{n/}{k}} \frac{ \eta_j( -v_t + (1-r) u_t - 3 \beta_0 r - s   )   }{A_{s, j}} e^{\frac n k h} \omega^n,\quad    \sup \psi_j = 0.$$
As $j\to \infty$ we have by the dominated convergence theorem
$$A_{s,j} =\frac { c_t^{\frac{n}{k}}}{V_t } \int_{X}\big( \eta_j( -v_t + (1-r ) u_t - 3 \beta_0 r - s   ) \big)  e^{\frac n k h} \omega^n \to 
A_s $$
where the constant $A_s$ is defined by
$$
A_s : = \frac { c_t^{\frac{n}{k}}}{V_t }\int_{\Omega_s} ( -v _t+ (1-r ) u_t - 3 \beta_0 r - s   ) e^{\frac n k h} \omega^n.
$$
Consider $\Phi = -\varepsilon \big(-\psi_j + \Lambda \big)^{\frac{n}{n+1}} + \big( - v _t+ (1-r) u_t - 3 \beta_0 r - s\big)$
where 
            $$\varepsilon =  \Big( \frac{k(n+1)}{n^2}   \Big)^{\frac{n}{n+1}} c(n,k)^{-\frac{1}{n+1}}   A_{s, j}^{\frac{1}{n+1 }},\quad  \mbox{where $c(n,k)$ is the one in (\ref{eqn:str})}$$
            $$ \Lambda = \Big( \frac{n}{n+1} \frac{\varepsilon}{r}  \Big)^{n+1} = C(n,k) \frac{A_{s, j}}{r^{n+1}} $$
Suppose $\sup\Phi = \Phi(x_0)$ for some point $x_0$ in $X$. If $x_0\not\in\Omega_s^\circ$, then by definition $\Phi(x_0)<0$. Otherwise $x_0\in\Omega_s^\circ$. We calculate as in \cite{GPT}. First note that for $G^{i\bar j} = \frac{\partial }{\partial (\omega_{t,v_t})_{i\bar j}} \log \sigma_k(\omega_{t, v_t})$
\begin{equation}\label{eqn:str}\det G^{i\bar j}\ge c(n, k) c_t^{-\frac n k } e^{-\frac n k h}
\end{equation} for some computable constant $c(n,k)>0$. It follows that at $x_0$ 
       \bea
       0& \ge &\nonumber  G^{i\bar j} (\Phi)_{i\bar j}(x_0)\\
       &\ge &\nonumber  \frac{n^2\varepsilon}{n+\delta} (-\psi_j + \Lambda)^{-\frac{1}{n+1}} \Big(\det G \cdot \det \omega_{t,\psi_j} \Big)^{1/n} - k + \Big( r - \frac{n\varepsilon}{n+1} \Lambda^{-\frac{1}{n+1}}   \Big) G^{i\bar j}( \omega_t)_{i\bar j}\\
         &\ge &\nonumber  \frac{n^2\varepsilon}{n+1} (-\psi _j+ \Lambda)^{-\frac{1}{n+1}}  c(n,k)^{1/n}  \Big(\frac{  -v _t+ (1-r) u_t - 3 \beta_0 r - s     }{A_{s,j}} \Big)^{1/n} - k.
        \eea
By the choice of $\varepsilon$ and $\Lambda$, we deduce that $\Phi(x_0)\le 0$. Thus $\Phi\le 0$ on $X$ and this implies 
\bea
\nonumber \int_{\Omega_s} \exp\Big\{ c_0  \Big(\frac{-v_t + (1-r) u_t - 3 \beta_0 r - s  }{A_{s, j}^{1/(n+1)}} \Big)^{\frac{n+1}{n}}     \Big\} \omega^n
&\leq&\int_{\Omega_s} \exp\Big(  - \alpha_0\psi_j + \alpha_0\Lambda    \Big)\omega^n
\nonumber\\
&\leq&  C \exp\Big\{\frac{A_{s, j}}{r^{n+1}}  \Big\}
\eea
for some small $c_0 = c_0(n,k,\omega, \chi)>0$, $C=C(n,k, \omega,\chi)>0$, and $\alpha_0$ a fixed number satisfying $0<\alpha_0<\alpha(X,\omega)$,
where $\alpha(X,\omega)$ is the $\alpha$-invariant of $(X,\omega)$. Letting $j\to\infty$ gives 
       \begin{equation}\label{eqn:useful}
        \int_{\Omega_s} \exp\Big\{ c_0  \Big(\frac{(-v _t+ (1-r) u_t - 3 \beta_0 r - s  )}{A_{s}^{1/(n+1)}} \Big)^{\frac{n+1}{n}}     \Big\} \omega^n
        \le C \exp\Big(\frac{A_{s}}{r^{n+1}} \Big).
      \end{equation}

\medskip

To estimate $\frac{A_s}{r^{n+1}}$ in (\ref{eqn:useful}), we need to consider separately the three cases {\bf I}, {\bf II} and {\bf III}.

\bigskip

\noindent $\bullet$ {\bf Case I}. In this case $k=n$ and $c_t = V_t$. By Lemma \ref{lemma Kol} we deduce that
       \bea
          \frac{A_{s}}{r^{n+1}} & = &\nonumber \frac{c_t^{n/k}}{V_t} \frac{1}{r^{n+1}} \int_{\Omega_s} ( -v_t + (1-r) u_t - 3 \beta_0 r - s  ) e^h    \omega^n\\
                                                 & \le &\nonumber  \frac{1}{r ^{n+1}} \int_{\Omega_s} ( -v _t+ (1-r) u_t - 3 \beta_0 r ) e^h   \omega^n\\
                                                 & \le  &\nonumber \frac{C(n,  \beta_0)}{r^{n+1}} \int_{E_2} e^h   \omega^n \le C_0 C(n,\beta_0) \gamma_1(r)r^{-1}\\
                                                 &\le  &\nonumber C(n,\beta_0) \quad \mbox{ by the choice of $\gamma_1(r)$.}
       \eea

 \medskip
 
 \noindent $\bullet$ {\bf Case II}. In this case, under our normalization, $V_1 = c_1 = 1$. As in Case {\bf I}, we have
        \bea
          \frac{A_{s}}{r^{n+1}} & =  &\nonumber \frac{1}{r^{n+1}} \int_{\Omega_s} ( -v_1 + (1-r) u_1 - 3 \beta_0 r - s  ) e^{\frac n kh}    \omega^n\\
                                                 & \le &\nonumber  \frac{1}{r ^{n+1}} \int_{\Omega_s} ( -v_1 + (1-r) u_1 - 3 \beta_0 r ) e^{\frac{n}{k} h}   \omega^n\\
                                                 & \le &\nonumber  \frac{C(n, \beta_0)}{r^{n+1}} \int_{E_2} e^{\frac {n}{k} h}   \omega^n =  \frac{C(n, \beta_0)}{r^{n+1}} \int_{E_2} e^{(\frac {n}{k} - 1) h}  e^h  \omega^n\\ %
                                                 &\le &\nonumber   \frac{C(n, \beta_0)}{r^{n+1}}\Big(  \int_{E_2} e^{ q_2^* (\frac {n}{k}-1) h} e^h   \omega^n     \Big)^{1/q_2^*} \Big(\int_{E_2} e^h   \Big)^{1/q_2}\\
                                                  &\le &\nonumber   \frac{C(n,  \beta_0)}{r^{n+1}}\Big(  \int_{E_2} e^{p_2 h}  \omega^n     \Big)^{1/q_2^*} \Big(\int_{E_2} e^h   \Big)^{1/q_2}\\
                                                  &\le &\nonumber  C(n,k, \beta_0, K) \gamma_2(r)^{\frac 1{q_2}} r^{\frac{n}{q_2} - n - 1} = C(n,k, \beta_0, K) , 
       \eea
where $\frac 1 {q_2} + \frac{1}{q_2^*} = 1$ and 
in the last equation we use the choice of the function $\gamma_2(r)$.

\bigskip

\noindent $\bullet$ {\bf Case III}. We note that since $[\chi]$ is big, $V_t \ge \int_X \chi^n >0$, hence $\frac{c_t^{n/k}}{V_t}\le C_{\omega, \chi}$ for a uniform $C_{\omega, \chi} = C_{\omega, \chi}(n, k)>0$ which we will fix throughout the proof below. Then we have
       \bea
          \frac{A_{s}}{r^{n+1}} & =& \nonumber \frac{c_t^{n/k}}{V_t} \frac{1}{r^{n+1}} \int_{\Omega_s} ( -v_t + (1-r) u_t - 3 \beta_0 r - s  ) e^{\frac nkh}    \omega^n\\
                                                       &\le& \nonumber   \frac{C(n, \omega,\chi, k, \beta_0)}{r^{n+1}}\Big(  \int_{E_2} e^{p_3 h}  \omega^n     \Big)^{1/q_3^*} \Big(\int_{E_2} e^h   \Big)^{1/q_3}\\
                                                  &\le& \nonumber C(n, k, \omega, \chi, K) \gamma_3(r)^{\frac 1{q_3}} r^{\frac{n}{q_3} - n - 1} = C(n, k, \omega, \chi, K), 
       \eea
where $\frac 1 {q_3} + \frac{1}{q_3^*} = 1$ and in the last identity we use the choice of the function $\gamma_3(r)$.

 \medskip

So for all cases {\bf I}, {\bf II} and {\bf III}, we get from (\eqref{eqn:useful}) that
       \begin{equation}\label{eqn:useful 1}
        \int_{\Omega_s} \exp\Big\{ c_0  \Big(\frac{-v _t+ (1-r) u_t - 3 \beta_0 r - s  }{A_{s}^{1/(n+1)}} \Big)^{\frac{n+1}{n}}     \Big\} \omega^n
        \le C , 
      \end{equation}
      for some constant $C>0$ depending on $n, k, \chi, \omega, K$ and the exponents $p_1, p_2, p_3$ in each case, respectively. In particular this $C$ is independent of the choice of $r\in (0,\bar r_0)$.

\bigskip
We choose $p>n$ as $p = p_1$ in case {\bf I}, and arbitrary and large $p>n$ in cases {\bf II} and {\bf III}.

\medskip 

Define $\eta: {{\mathbb R}}_+\to {{\mathbb R}}_+$ by $\eta(x) = (\log(1+x))^p$. Note that $\eta$ is a strictly increasing function with $\eta(0) = 0$, and let $\eta^{-1}$ be its inverse function. If we let
\bea
\Psi: = \frac{c_0}{2}\Big(\frac{-v _t+ (1-r) u_t - 3 \beta_0 r - s  }{A_{s}^{1/(n+1)}} \Big)^{\frac{n+1}{n}} 
\eea
then we have for any $z\in\Omega_s$, by
the generalized Young's inequality with respect to $\eta$, 
   \bea
  \Psi(z)^p e^{\frac{n}{k} h(z)} & \le & \nonumber \int_0^{\exp({\frac{n}{k} h(z)})} \eta(x) dx + \int_0^{\Psi(z)^p} \eta^{-1}(y) dy\\
    &\le & \nonumber e^{\frac{n}{k} h(z)} (1+  |h(z)|)^p + C(p) e^{2\Psi(z)}
    \eea
We integrate both sides in the inequality above over $z\in \Omega_s$, and get by (\eqref{eqn:useful 1}) that 
   \bea
  \int_{\Omega_s}  \Psi(z)^p e^{\frac{n}{k} h(z)} \omega^n  
   \le\nonumber \| e^{h} \|_{L^{n/k}(\log L)^p} + C,
  \eea
  where the constant $C>0$ depends only on $n,k,\omega, \chi, p, K$. In view of the definition of $\Psi$, this implies  
           \begin{equation}\label{eqn:final 1}
     \int_{\Omega_s} (-v _t+ (1-r) u_t - 3 \beta_0 r - s  )^{\frac{(n+1)p}{n}} e^{\frac{n}{k}h } \omega^n 
   \le  C A_s^{\frac{p}{n}} \big(\| e^{h} \|_{L^{n/k}(\log L)^p} + 1\big).
             \end{equation}
             %
It follows from the H\"older inequality that
\bea
     A_s  & = &\nonumber  \nonumber\frac{c_t^{n/k}}{V_t} \int_{\Omega_s} (-v _t+ (1-r) u_t - 3 \beta_0 r - s  ) e^{\frac n k h} \omega_X^n \\
     & \le\nonumber  &\nonumber  \Big( \frac{c_t^{n/k}}{V_t} \int_{\Omega_s} ( -v _t+ (1-r) u_t - 3 \beta_0 r - s   ) ^{\frac{(n+1)p}{n}} e^{\frac n k h}\omega^n \Big)^{\frac{n}{(n+1)p}} \cdot \Big(\frac{c_t^{n/k}}{V_t}  \int_{\Omega_s} e^{\frac n k h} \omega ^n \Big)^{1/q}\\
     &\le &\nonumber C  A_s^{\frac {1}{n+1}}\Big( \| e^{h} \|_{L^{n/k}(\log L)^p} + 1)\Big)^{\frac{n}{(n+1)p}} \cdot \Big(\frac{c_t^{n/k}}{V_t}  \int_{\Omega_s} e^{\frac n k h} \omega^n \Big)^{1/q}
       \eea
where $q>1$ satisfies  $\frac{n}{p(n+1)} + \frac{1}{q} = 1$, i.e. $q = \frac{p(n+1)}{p(n+1) - n}$. The inequality above yields 
\begin{equation}\label{eqn:1.21}
   A_s\le C \Big(\| e^{h} \|_{L^{n/k}(\log L)^p} + 1 \Big)^{1/p} \cdot \Big(\frac{c_t^{n/k}}{V_t}  \int_{\Omega_s} e^{\frac n k h} \omega^n \Big)^{\frac{1+n}{qn}} = B_0 \Big(\frac{c_t^{n/k}}{V_t}  \int_{\Omega_s} e^{\frac n k h} \omega^n \Big)^{1+\delta_0}.
\end{equation}
Observe that the exponent of the integral on the right hand of (\eqref{eqn:1.21}) satisfies
$$\frac{1+n}{qn} = \frac{pn + p - n}{pn} = 1+ \delta_0>1, \quad \mbox{for $\delta_0: = \frac{p-n}{pn}>0$.}$$
We remark that $\delta_0$ can be chosen to be close to $1/n$ in cases {\bf II} and {\bf III} by picking $p$ large enough. 
 Furthermore, we note that \begin{equation}\label{eqn:B0}B_0 := C \Big(\| e^{h} \|_{L^{n/k}(\log L)^p} + 1 \Big)^{1/p}\end{equation} is a constant depending only on $n,k, \omega, \chi, K$, and the exponents $p_1,p_2,p_3$ in each case, respectively, and in particular, it is independent of $r$ with $r\in (0,\bar r_0)$. 

     If we define
     $$\phi(s) =\frac{c_t^{n/k}}{V_t} \int_{\Omega_s} e^{\frac n k h} \omega^n,$$  then (\eqref{eqn:1.21}) shows that if $\Omega_{s+s'}\not = \emptyset$ then
     \begin{equation}\label{eqn:iteration}
     s' \phi(s+s')\le  B_0 \phi(s)^{1+\delta_0},\quad \mbox{for all } s' \ge 0 \mbox{ and }s\ge 0.\end{equation}
    
We now choose $r_0< \bar r_0$ small in each case as follows.
 
 \medskip
 \noindent{\bf Case I}. We choose $r_0>0$ small so that for $r\in (0,r_0)$
      \bea
  B_0 \phi(0)^{\delta_0}\le B_0\Big (\int_{E_2} e^{ h} \omega^n\Big)^{\delta_0}\le B_0 C_0^{\delta_0}  (\gamma_1(r) r^n)^{\delta_0} \le B_0 C_0^{\delta_0}  (\gamma_1(r_0) r_0^n)^{\delta_0}\le \frac 1 2
        \eea
  and $\phi(0)\le C_0 \gamma_1(r) r^n< \bar  C r^{1/\delta_0}$ by Lemma \ref{lemma Kol} for some uniform $\bar C$. 
  
  \medskip
  
  \noindent{\bf Case II}. We choose $r_0>0$  small so that for all $r\in (0,r_0)$ 
      \bea
  B_0 \phi(0)^{\delta_0}& \le &\nonumber  B_0\Big (\int_{E_2} e^{ \frac n kh} \omega^n\Big)^{\delta_0} \\ 
  & \le &\nonumber  B_0 \Big(  \int_{E_2} e^{p_2 h}  \omega^n     \Big)^{\delta_0 /q_2^*} \Big(\int_{E_2} e^h  \omega^n \Big)^{\delta_0/q_2}\\
 &  \le &\nonumber  B_0 C_0^{\delta_0/q_2 } \Big(  \int_{E_2} e^{p_2 h}  \omega^n     \Big)^{\delta_0 /q_2^*} (\gamma_2(r) r^n)^{\delta_0/q_2}  \\
  &  \le &\nonumber  B_0 C_0^{\delta_0/q_2 } \Big(  \int_{E_2} e^{p_2 h}  \omega^n     \Big)^{\delta_0 /q_2^*} (\gamma_2(r_0) r_0^n)^{\delta_0/q_2} \le  \frac  12 
        \eea
   where $\frac 1 q _2+ \frac{1}{q^*_2} = 1$ and we also have 
   $$\phi(0)\le C_0^{1/q_2} \Big(  \int_{E_2} e^{p_2 h}  \omega^n     \Big)^{1 /q_2^*} (\gamma_2(r) r^n)^{1/q_2} < \bar C r^{1/\delta_0}$$
  for some uniform $\bar C>0$ by the definition of $\gamma_2(r)$.
    
  \medskip
  
  \noindent{\bf Case III}. We choose $r_0>0$ small so that for $r\in (0,r_0)$
      \bea
  B_0 \phi(0)^{\delta_0}& \le &\nonumber  B_0 C_{\omega,\chi}^{\delta_0}\Big (\int_{E_2} e^{ \frac n kh} \omega^n\Big)^{\delta_0} \\ 
  & \le &\nonumber  B_0C_{\omega,\chi}^{\delta_0} \Big(  \int_{E_2} e^{p_3 h}  \omega^n     \Big)^{\delta_0 /q_3^*} \Big(\int_{E_2} e^h  \omega^n \Big)^{\delta_0/q_3}\\
 &  \le  &\nonumber B_0C_{\omega,\chi}^{\delta_0} C_0^{\delta_0/q_3 } \Big(  \int_{E_2} e^{p_3 h}  \omega^n     \Big)^{\delta_0 /q_3^*} (\gamma_3(r) r^n)^{\delta_0/q_3}  \\
  &  \le &\nonumber  B_0C_{\omega,\chi}^{\delta_0} C_0^{\delta_0/q_3 } \Big(  \int_{E_2} e^{p_3 h}  \omega^n     \Big)^{\delta_0 /q_3^*} (\gamma_3(r_0) r_0^n)^{\delta_0/q_3}\le  \frac  12 
        \eea
   where $\frac 1 q _3+ \frac{1}{q^*_3} = 1$ and we also have 
   $$\phi(0)\le C_{\omega,\chi} C_0^{1/q_3} \Big(  \int_{E_2} e^{p_3 h}  \omega^n     \Big)^{1 /q_3^*} (\gamma_3(r) r^n)^{1/q_3} < \bar C r^{1/\delta_0}$$
   by the choice of $\gamma_3(r)$.
  
\medskip

It is clear that in all cases, $r_0$ and $\bar C$ depend only on the given data, namely, $n, k, \omega,\chi, K$ and $p_a$, and we have $B_0 \phi(0)^{\delta_0} \le \frac 1 2$ and $\phi(0)\le \bar C r^{1/\delta_0}$. 

\medskip

 Define a sequence of increasing real numbers $(s_j)$ inductively such that  $s_0 = 0$ and 
 $$s_{j+1} = \sup\{s>s_j| \phi(s)> \frac 1 2 \phi(s_j)\}.$$
 Then we can show  that (see \cite{GPT})  $\phi(s_j)\le 2^{-j} \phi(s_0)$ and 
 $$s_{j+1} - s_j \le 2 B_0 2^{-j \delta_0} \phi(0)^{\delta_0}. $$
 
 Thus the limit $S_\infty =\lim_{j\to\infty}s_j$ satisfies 
 $$S_\infty \le \frac{2B_0}{1-2^{-\delta_0}} \phi(0)^{\delta_0}\le \frac{2B_0 \bar C^{\delta_0}}{1-2^{-\delta_0}} r = \hat C r.$$ %
 
 Hence the set $ \Omega_{\hat C r} = \emptyset$, and we conclude that 
 $$v_t\ge (1-r) u_t - 3 \beta_0 r -\hat Cr, \quad \mbox{or equivalently }\quad v_t - u_t \ge - Cr,$$ for some uniform constant $C>0$ depending only on the given data.
 By the normalization $\max(u_t-v_t) = \max(v_t-u_t)$, it is clear that $v_t-u_t\le Cr$. The proof of Theorem \ref{thm:main} is complete.


\smallskip


\bigskip


\noindent Department of Mathematics \& Computer Science, Rutgers University, Newark, NJ 07102 USA

\noindent bguo@rutgers.edu,

\medskip

\noindent Department of Mathematics, Columbia University, New York, NY 10027 USA

\noindent phong@math.columbia.edu, tong@math.columbia.edu

\end{document}